\documentclass[12pt]{article}
\usepackage{mathrsfs}
\usepackage{amsfonts}
\usepackage{latexsym}
\usepackage{amscd}
\date{}

\begin{document}
\title{  Finitistic dimension and Endomorphism algebras of Gorenstein projective modules}
\author{{\small Aiping Zhang$^{1,2}$}\\
{\small $^1$School of Mathematics,\ Shandong University,\ Jinan
250100,\ China }\\
{$^2$\small School of Mathematics and Statistics,\ Shandong University,\ Weihai,\
Weihai 264209,\ China }\\
\\
}

\pagenumbering{arabic}

\maketitle

\footnote{ {\it Email addresses}: aipingzhang@sdu.edu.cn}
\begin{center}
 \begin{minipage}{120mm}
   \small\rm
  Let $A$ be an Artin algebra, $M$ be a Gorenstein projective $A$-module and $B =$ End$_A M$, then $M$ is a $A$-$B$-bimodule. We use the restricted flat dimension of $M_B$ to give a characterization of the homological dimensions of $A$ and $B$, and obtain the following main results: (1) if $A$ is  a CM-finite algebra with  $\cal GP$($A$) = add$_AE$ and fin.dim $A \geq 2,$  then ${\rm fin.dim}\ B \leq {\rm fin.dim}\ A + {\rm rfd}(M_B) +{\rm pd}_B{\rm Hom}_A(M, E);
  $ (2) If $A$ is  a CM-finite $n$-Gorenstein algebra with $\cal GP$($A$) = {\rm add}$_AE$ and $n \geq 2$, then {\rm gl.dim }$B \leq n + {\rm pd}_B{\rm Hom}_A(M, E).$
\end{minipage}
\end{center}
\begin{center}
  \begin{minipage}{120mm}
\end{minipage}
\end{center}
\vskip -0.3cm
\begin{center}
\begin{minipage}{120mm}
\end{minipage}
\end{center}

\section {Introduction}

\vskip 0.2in

\ \ \ \ \ \

In 1969, M.Auslander and M. Bridger [1] introduced the modules of $G$-dimension zero over two-sided noetherian rings, which can be viewed as a generalization of finitely generated projective modules. Several decades later, Enochs and Jenda [7, 8, 10] extended the ideas
of Auslander and Bridger and introduced the notion of Gorenstein projective modules (not necessarily finitely generated) over a general ring. It is known that these two notions coincide
for finitely generated modules over a left and right noetherian ring. The class of Gorenstein projective modules has become a main ingredient in the relative homological algebra, and is widely used in the representation theory and algebraic geometry [3-5, 11, 12].

Let $A$ be an Artin algebra. The famous finitistic dimension conjecture says that there exists a uniform bound of the finite projective dimensions of all finitely generated left $A$-modules of finite projective dimension. This conjecture is related to many other homological conjectures and attracts many algebraists. Motivated by a result of Auslander, which states that every Artin algebra $B$ is of the form $eAe$ with $A$ being of finite global dimension, there is much literature on the study of the relationship of homological dimensions of $A$ and $eAe$. Xi [19] showed under the condition the Gorenstein projective dimension of right $eAe$-module $Ae$ is finite, if fin.dim $A$ is finite, then fin.dim $eAe$ is finite. Wei [18] used the restricted flat dimension of right $eAe$-module $Ae$ to characterize  the finitistic dimensions of $A$ and $eAe$. Huang and Huang [14] used the projective dimension of right $eAe$-module $Ae$ to characterize  the $m$-syzygy finiteness of $A$ and $eAe$. For more details, please see [14,15,17-20 ] and the references therein.

In [16], Li and Zhang showed that if $A$ is a CM-finite $n$-Gorenstein algebra with $\cal GP$($A$) = {\rm add}$_AE$ and $n \geq 2$, then the global dimension of End$_AE$ is no more than $n$. However, in general case, we do not know much about the relationship of the homological dimensions of $A$ and the endomorphism algebras of Gorenstein projective $A$-modules. Motivated by the above statements, the aim of this paper is to give a characterization  of homological dimensions of $A$ and the endomorphism algebra $B$ of a Gorenstein projective $A$-module $M$. Our main results are as follows.

 \vskip 0.2in

 {\bf Theorem 1.1.}  {\it  Let $M$ be a Gorenstein projective $A$-module and $B = {\rm End}_AM$. If $A$ is  CM-finite with $\cal GP$($A$) = {\rm add}$_AE$ and ${\rm fin.dim}\ A \geq 2$, Then ${\rm fin.dim}\ B \leq {\rm fin.dim}\ A + {\rm rfd}(M_B) + {\rm pd}_B{\rm Hom}_A(M, E).$}

 \vskip 0.2in

 {\bf Theorem 1.2.}  {\it Let $M$ be a Gorenstein projective $A$-module and $B = {\rm End}_AM$. If $A$ is  a CM-finite $n$-Gorenstein algebra with $\cal GP$($A$) = {\rm add}$_AE$ and $n \geq 2$, then {\rm gl.dim }$B \leq n + {\rm pd}_B{\rm Hom}_A(M, E).$}

 \vskip 0.2in

This paper is organized as follows: In section 2, we recall some definitions and basic facts. The proofs of the results will be given in section 3.¡¡

 \vskip 0.2in

\section { Preliminaries}

 \vskip 0.2in

\ \ \ \ \ \

Let $A$ be a ring, and $T_A$ be a right $A$-module. Following by [6], $T_A$ is said to have (big) restricted flat dimension at most $m$ if for each
$i > m$ the functor Tor$^A_{i}(T, -)$ vanishes on the category of modules of finite flat dimension. The little restricted flat dimension is defined correspondingly by considering only modules of finite flat dimension which admit a projective resolution of with finitely generated projectives. The restricted flat (resp., little restricted flat) dimension of $T_A$ is denoted by Rfd$T_A$(resp., rfd$T_A$).

Throughout this paper, let $A$ be an Artin algebra. We denote by $A$-mod the category of all finitely generated left $A$-modules and $\cal P$$(A)
$ the subcategory of projective $A$-modules. Given an $A$-module $M$, we denote by add$_AM$ the
category of all direct summands of finite direct sums of $_AM$. The nth syzygy of the $A$-module $M$ is denoted by $\Omega_A^n(M)$. Let $\cal C$ be a subcategory of $A$-mod and $M\in A$-mod, we denote by $\cal C$-dim$(_AM)$ the minimal integer $m$ such that there is an exact sequence $0 \rightarrow T_m \rightarrow \cdots \rightarrow T_0 \rightarrow M \rightarrow 0$ with each $T_i\in \cal C$ and call it the $\cal C$-dimension of $_AM$.

Let $A$ be an Artin algebra. An $A$-module $M$ is said to be $Gorenstein$ $projective$ if there is an exact sequence of projective modules in $A$-mod,
$$ \cdots \rightarrow P_1 \rightarrow P_0 \rightarrow P^0 \rightarrow P^1 \rightarrow  \cdots, $$
such that $M \cong {\rm Im} (P_0 \rightarrow P^0)$ and Hom $(- ,\ Q)$ is still exact for every projective module $Q \in A$-mod.

We say that an $A$-module $M\in A$-mod has $Gorenstein$-$dimension$ $n$ if there is an exact sequence
$0 \rightarrow X_n \rightarrow \cdots \rightarrow X_1 \rightarrow X_0 \rightarrow M \rightarrow 0$ such that all $X_i$ are Gorenstein projective and that $n$ is minimal among the lengths of such exact sequence. In this case we write Gpd$_AM = n$.

A Gorenstein ring $R$ is a left and right Noetherian ring with id $_RR < \infty$ and id $R_R < \infty$; a Gorenstein ring $R$ is $n$-Gorenstein if id $_RR \leq n$, and in
this case id $R_R \leq n$. An algebra $A$ is Gorenstein if it is a Gorenstein ring.

Let $\cal GP$($A$) be the full subcategory of $A$-mod of Gorenstein projective modules. An algebra $A$ is called CM-$finite$, if $\cal GP$($A$) has only finitely
many isomorphism classes  of indecomposable objects. Clearly, $A$ is CM-$finite$ if and only if there is an $A$-module $E$ such that $\cal GP$($A$) = add$_AE$. If gl.dim $A < \infty$, then $\cal GP$($A$) = $\cal P$($A$).

The following lemma proved in [21] will be used later.

     \vskip 0.2in

     {\bf Lemma 2.1.}  {\it Let $M$ be a Gorenstein projective $A$-module and $B = {\rm End}_AM$. Then for any $X\in B$-{\rm mod}, $\Omega_B^2(X) \cong {\rm Hom}_A(M, \Omega_A^2(N)\oplus Q)$ for some $N \in A$-{\rm mod} with {\rm Gpd}$_A(M\otimes_BX) = {\rm Gpd}(_AN)$ and $_AQ$ is a projective $A$-module depending on $X$ }.

\vskip 0.2in

Lemma 2.2. is known in homological algebra.

 {\bf Lemma 2.2.}  {\it Let $A$ be an Artin algebra and let $M$ be an $A$-module. If there is an exact sequence
  $$0 \rightarrow X_s \rightarrow  \cdots \rightarrow X_1 \rightarrow X_0 \rightarrow M \rightarrow 0$$
  of $A$-modules with {\rm pd}$_AX_i \leq k$ for all $i$, then {\rm pd}$_AM \leq s + k$.}

\vskip 0.2in

In this paper, we follow the standard terminology and notation used in the representation theory of algebras and relative
homological algebra, see [2, 9].

\section { Finitistic dimensions and Endomorphism algebras of Gorenstein projective modules}

\vskip 0.2in

\ \ \ \ \ \

Let $M$ be an $A$-module with $B = $End $_AM$. Then $M$ is also a right $B$-module. In this section, We use the restricted flat dimension of $M_B$ to give a characterization of the homological dimension of $A$ and $B$.

\vskip 0.2in

{\bf Theorem 3.1.}  {\it  Let $M$ be a Gorenstein projective $A$-module and $B = {\rm End}_AM$. If $A$ is  CM-finite with $\cal GP$($A$) = {\rm add}$_AE$ and ${\rm fin.dim}\ A \geq 2$,, Then ${\rm fin.dim}\ B \leq {\rm fin.dim}\ A + {\rm rfd}(M_B) + {\rm pd}_B{\rm Hom}_A(M, E).$}

\vskip 1mm

{\bf Proof.} If fin.dim $A$ or rfd$(M_B)$ is infinite, then we have nothing to say. So we assume that fin.dim $A = r < \infty$ and rfd$(M_B) = t < \infty$.

Let $_BY \in B$-mod with pd$(_BY) < \infty$. Denote by $Y_i$ the ith syzygy of $Y$, for each i.
Since rfd($M_B) = t$, we have Tor$^B_{i+t}(M, Y)$ = Tor$^B_{i}(M, Y_t) = 0$ for $i \geq 1$.

Let $0 \rightarrow P_m\rightarrow \ldots\rightarrow P_1\rightarrow P_0\rightarrow Y_t\rightarrow 0$ be a projective resolution of $Y_t$ in $B$-mod.
By applying the functor $M\otimes_B-$ to the
above exact sequence and Tor$^B_{i}(M, Y_t) = 0$ for $i \geq 1$, we obtain the following exact sequence

$$ 0 \rightarrow M\otimes_BP_m\rightarrow \ldots\rightarrow M\otimes_BP_1\rightarrow M\otimes_BP_0\rightarrow M\otimes_BY_t\rightarrow 0.\ \ (*)$$

Now consider the exact sequence $(*)$, since $M\otimes_BP_i(i = 0, 1, \ldots, m) \in$ add$_AM$ and the subcategory of Gorenstein projective modules is closed under direct sums and direct summands, we conclude that $ M\otimes_BP_i$ are still Gorenstein projective and thus Gpd$_A(M\otimes_BY_t) < \infty$.
It follows from [12, Theorem 2.28] that Gpd$_A(M\otimes_BY_t) \leq r$.

By Lemma 2.1, $\Omega_B^2(Y_t) \cong {\rm Hom}_A(M, \Omega_A^2(N)\oplus Q)$ for some $N \in A$-{\rm mod} with {\rm Gpd}$_A(M\otimes_BY_t) = {\rm Gpd}_A(N)$ and $_AQ$ is a projective $A$-module depending on $Y_t$. Since Gpd$_A(M\otimes_BY_t) \leq r$, we obtain that Gpd$_A(\Omega_A^2(N)) \leq r-2.$

By [13, Theorem 2.6], there exists an exact sequence
$$0\rightarrow Q_{r-2} \rightarrow \ldots \rightarrow Q_1 \rightarrow G_0 \rightarrow \Omega_A^2(N)\oplus Q \rightarrow 0\ \  (**)$$
in $A$-mod where $G_0$ is Gorenstein projective and $Q_i$ is projective for $i=1, \ldots, r-2.$ Moreover by applying the functor Hom$_A(M,-)$ to $(**)$, we obtain the
following exact sequence:

$ 0\rightarrow {\rm Hom}_A(M, Q_{r-2})\rightarrow \ldots \rightarrow {\rm Hom}_A(M, Q_1)\rightarrow {\rm Hom}_A(M, G_0)\rightarrow {\rm Hom}_A(M,\\ \Omega_A^2(N)\oplus Q ) \rightarrow 0.$
Thus we have

$$
\begin{array}{ccl}
{\rm pd}\ _{B}Y&\leq & {\rm pd}\ \Omega_B^2(Y_t) + t + 2\\
     &= & {\rm pd}_B{\rm Hom}_A(M, \Omega_A^2(N)\oplus Q ) + t + 2\\
     &\leq & r + t + {\rm pd}_B{\rm Hom}_A(M, E).\\
    \end{array}
     $$
The proof is completed.
\hfill$\Box$

\vskip 0.2in
Note that in Theorem 3.1., if fin.dim $A = 0$ or $1$, then  ${\rm fin.dim}\ B \leq 2 + {\rm rfd}(M_B) + {\rm pd}_B{\rm Hom}_A(M, E).$

\vskip 0.2in

{\bf Theorem 3.2.}  {\it Let $M$ be a Gorenstein projective $A$-module and $B = {\rm End}_AM$. If $A$ is  a CM-finite $n$-Gorenstein algebra with $\cal GP$($A$) = {\rm add}$_AE$ and $n \geq 2$, then {\rm gl.dim }$B \leq n + {\rm pd}_B{\rm Hom}_A(M, E).$}

\vskip 1mm

{\bf Proof.}
Let $_BY \in B$-mod. By Lemma 2.1, $\Omega_B^2(Y) \cong {\rm Hom}_A(M, \Omega_A^2(N)\oplus Q)$ for some $N \in A$-{\rm mod} with {\rm Gpd}$_A(M\otimes_BY) = {\rm Gpd}_A(N)$ and $_AQ$ is a projective $A$-module depending on $Y$. Since $A$ is a CM-finite $n$-Gorenstein algebra,  we have Gpd$_A(\Omega_A^2(N)) \leq n-2.$

By [13, Theorem 2.6], there exists an exact sequence
$$0\rightarrow Q_{n-2} \rightarrow \ldots \rightarrow Q_1 \rightarrow G_0 \rightarrow \Omega_A^2(N)\oplus Q \rightarrow 0\ \  (1)$$
in $A$-mod where $G_0$ is Gorenstein projective and $Q_i$ is projective for $i=1, \ldots, n-2.$ Moreover by applying the functor Hom$_A(M,-)$ to the exact sequence $(1)$, we obtain the
following exact sequence:

$ 0\rightarrow {\rm Hom}_A(M, Q_{n-2})\rightarrow \ldots \rightarrow {\rm Hom}_A(M, Q_1)\rightarrow {\rm Hom}_A(M, G_0)\rightarrow {\rm Hom}_A(M,\\ \Omega_A^2(N)\oplus Q ) \rightarrow 0.$
Thus we have

$$
\begin{array}{ccl}
{\rm pd}\ _{B}Y&\leq & {\rm pd}\ \Omega_B^2(Y) + 2\\
     &= & {\rm pd}_B{\rm Hom}_A(M, \Omega_A^2(N)\oplus Q ) + 2\\
    &\leq & n + {\rm pd}_B{\rm Hom}_A(M, E).\\
    \end{array}
     $$
 The proof is completed.
\hfill$\Box$

\vskip 0.2in

Note that in Theorem 3.2., if $n = 0$ or $1$, then gl.dim $B \leq 2 + {\rm pd}_B{\rm Hom}_A(M, E)$.

\vskip 0.2in

As applications of Theorem 3.2, we can easily obtain the following two corollaries. Note that Corollary 3.1. is also a result of [16].

\vskip 0.2in

{\bf Corollary 3.1.}  {\it Let $A$ be a CM-finite $n$-Gorenstein algebra with $\cal GP$($A$) = {\rm add}$_AE$ and $n \geq 2$, then {\rm gl.dim End}$_AE \leq n$}.

\vskip 0.2in

{\bf Corollary 3.2.}  {\it Let $A$ be an Artin algebra and $e$ be an idempotent of $A$. If ${\rm gl.dim}\ A \geq 2$, then \it {\rm gl.dim} $eAe \leq$ {\rm gl.dim} $A +$ {\rm pd}$_{eAe} eA$}.

\vskip 0.2in

The following result is also proven in [18], here we give another proof. Note that Theorem 3.3. contains [19, Theorem 1.2].

{\bf Theorem 3.3.}  {\it Let $A$ be an Artin algebra and $e$ be an idempotent of $A$. If ${\rm fin.dim}\ A \geq 2$, then {\rm fin.dim} $eAe \leq$ {\rm fin.dim} $A + ${\rm rfd}$(Ae_{eAe})$}.

\vskip 1mm

{\bf Proof.} We denote by fin.dim $A = n$ and $B = eAe$. The first part of the proof is similar to that of Theorem 3.1. For completeness, we repeat it here.

Let $_BY \in B$-mod with pd$(_BY) < \infty$. Denote by $Y_i$ the ith syzygy of $Y$, for each i.
Since rfd($Ae_B) = t$, we have Tor$^B_{i+t}(Ae, Y)$ = Tor$^B_{i}(Ae, Y_t) = 0$ for $i \geq 1$.

Let $0 \rightarrow P_m\rightarrow \ldots\rightarrow P_1\stackrel{f_1}\rightarrow P_0\stackrel{f_0} \rightarrow Y_t\rightarrow 0$ be a projective resolution of $Y_t$ in $B$-mod.
By applying the functor $Ae\otimes_B-$ to the
above exact sequence and Tor$^B_{i}(Ae, Y_t) = 0$ for $i \geq 1$, we obtain the following exact sequence

$$ 0 \rightarrow Ae\otimes_BP_m \rightarrow  \ldots \rightarrow  Ae\otimes_BP_1 \stackrel{1\otimes f_1}\rightarrow Ae\otimes_BP_0 \stackrel{1\otimes f_0} \rightarrow Ae\otimes_BY_t\rightarrow 0.\ \ (2)$$

We denote by $K = {\rm ker}(1\otimes f_1)$, then through the proof of Lemma 2.1., we know $\Omega_B^2(Y_t) \cong {\rm Hom}_A(Ae, K)$. Since $Ae\otimes_BP_i \in$ add$_AAe$ for $i = 0,\cdots,m$, then from exact sequence (2), we conclude pd$_AK$ is finite, pd$_AK \leq n-2$ and add$Ae$-dim($_AK$) is finite. Hence there exists an exact sequence

$$ 0 \rightarrow Q_{n-2} \rightarrow \cdots \rightarrow Q_1 \rightarrow Q_0 \rightarrow K \rightarrow 0$$

in $A$-mod such that $Q_i\in {\rm add}_AAe$ for $i = 0, 1, \cdots, n-2$. Moreover by applying the functor Hom$_A(Ae,-)$ to the above exact sequence, we obtain the following exact sequence

$ 0\rightarrow {\rm Hom}_A(Ae, Q_{n-2})\rightarrow \ldots \rightarrow {\rm Hom}_A(Ae, Q_1)\rightarrow {\rm Hom}_A(Ae, Q_0)\rightarrow {\rm Hom}_A(Ae,K)\rightarrow 0.$
Thus we have

$$
\begin{array}{ccl}
{\rm pd}\ _{B}Y&\leq & {\rm pd}\ \Omega_B^2(Y_t) + t + 2\\
     &= & {\rm pd}_B{\rm Hom}_A(Ae,K ) + t + 2\\
     &\leq & n + t\\
    \end{array}
     $$
since ${\rm Hom}_A(Ae, Q_i)$ for $i = 0,\cdots,n-2$ are all projective $B$-modules. The proof is completed.
\hfill$\Box$

\vskip 0.2in

\section {Acknowledgements}

 The research was supported by the National Nature Science Foundation of China (Grant No.11601274).

\begin{description}

\item{[1]}\ M.Auslander, M.Bridger, Stable Module theory, Mem. Amer. Math. Soc., vol.94, Amer. Math. Soc., Providence, RI, 1969.

\item{[2]}\ M.Auslander, I.Reiten, S.O.Smal$\phi$,  Representation Theory of Artin Algebras, Cambridge Univ. Press, 1995.

\item{[3]}\ L.L.Avramov, A. Martsinkovsky, Absolute, relative, and Tate cohomology of modules of finite Gorenstein dimension, Proc. London Math.Soc.(3) 85 (2002) 393-440.

\item{[4]}\ L.W.Christensen, Gorenstein Dimension, Lecture Notes in Math., vol.1747, Springer-Verlag, Berlin, 2000.

\item{[5]}\ L.W.Christensen, A.Frankild, H.Holm, On Gorenstein projective, injective and flat dimensions-a functorial description with applications, J.Algebra. 302 (2006) 231-279.

\item{[6]}\ L.W.Christensen, H.-B. Foxby, A.Frankild, Restricted homological dimensions and Cohen-Macaulayness, J.Algebra 251(1)(2002)479-502.

\item{[7]}\ E.E.Enochs, O.M.G.Jenda, On Gorenstein injective modules, Comm.Algebra. 21 (10) (1993) 3489-3501.

\item{[8]}\ E.E.Enochs, O.M.G.Jenda, Gorenstein injective and projective modules, Math.Z. 220 (1995) 611-633.

\item{[9]}\ E.E.Enochs, O.M.G.Jenda, Relative Homological Algebra, de Gruyter Exp.Math., vol.30, Walter de Gruyter, Berlin, New York, 2000.

\item{[10]}\ E.E.Enochs, O.M.G.Jenda, J.Xu, Foxby duality and Gorenstein injective and projective modules, Trans.Amer.Math.Soc. 348 (8) (1996) 3223-3234.

\item{[11]}\ H.-B.Foxby, Gorenstein dimension over Cohen-Macaulay ring, in :W.Bruns(Ed.), Proceedings of International Conference on Commutative Algebra,Universit\"{a}t Onsabr\"{u}ck, 1994.

\item{[12]}\ H.Holm, Gorenstein homological dimensions, J.Pure Appl.Algebra. 189 (2004) 167-193.

\item{[13]}\ C.Huang, Z.Huang, Gorenstein syzygy modules, J.Algebra. 324 (2010) 3408-3419.

\item{[14]}\ Z.Huang, J.Sun, Endomorphism algebras and Igusa-Todorov algebras. Acta Math.Hungar. 140(1-2)(2013)60-70.

\item{[15]}\ K.Igusa, G.Todorov, On the finitistic global dimension conjecture for Artin algebras, in: $Representations\ of\ Algebras\ and\ Related\ Topics$, Fields Inst.Commun.45, Amer.Math.Soc.(Province, RI, 2005), pp. 201-204.

 \item{[16]}\ Z.Li, P.Zhang, Gorenstein algebras of finite Cohen-Macaulay type, Adv. Math. 223(2010)728-734.

\item{[17]}\ J. Wei, Finitistic dimensionn and Igusa-Todorov algebras, Adv. Math. 222 (2009) 2215-2226.

\item{[18]}\ J. Wei, Finitistic dimension and restricted flat dimension, J. Algebra. 320(2008)116-127

\item{[19]}\ C.Xi, On the finitistic dimension conjecture III: Related to the pair $eAe \subseteq A$, J.Algebra. 319 (2008) 3666-3668.

\item{[20]}\ A.Zhang, S.Zhang, on the finitistic dimension conjecture of artin algebras, J.Algebra. 320 (2008) 253-258.

\item{[21]}\ A.Zhang, Endomorphism algebras of Gorenstein projective modules, J.Algebra Appl.17(9)(2018)1850177.

\

\end{description}

\end{document}